# ONLINE RESOURCES IN MATHEMATICS: TEACHERS' GENESIS OF USE

Laetitia Bueno-Ravel, Ghislaine Gueudet

CREAD, IUFM Bretagne, France

*The long-term objective of our research is to develop the instrumental approach for teachers. A first step, presented in this paper, is to observe stable behaviours of teachers using internet resources in mathematics. We retain the scenarios as indicators of the genesis processes. We propose a scenario taxonomy taken from categories elaborated by computer sciences specialists and complemented to take into account didactical aspects. The descriptions provided by teachers permitted to observe an evolution of their scenarios elaboration's practices.*

## INTRODUCTION

We present here the work of a team comprising two mathematics education researchers, two teacher trainers and five teachers, and belonging to a more general research project termed GUPTEN (which holds for, in French: Genesis of Professional Uses of Technology by Teachers, project directed by Jean-Baptiste Lagrange).

This project aims at studying, in particular, the way teachers take over a new technological tool. The word "Genesis" in GUPTEN stems from the instrumental approach, a theoretical frame encompassing elements both from cognitive ergonomy and from the anthropological theory of didactics (Chevallard 1999) and developed within mathematics didactics to study issues related with ICT. Rabardel (1995) stresses the difference between an artifact, which is a given object, and an instrument. The instrument is a psychological construct; it comprises the artifact, and the schemes and techniques the user develops while using it (Guin and Trouche 2002).The building process of the instrument is called "the instrumental genesis".

We study here specific technological resources: e-exercises bases (shortened EEB in what follows). These resources consist of exercises classified according to their mathematical content, to their difficulty, and/or to the mathematical tools they require. These exercises are associated with an environment which consists of suggestions, correction, explanations, tools for the resolution of the exercise, score etc. (for more details about the possible features of an e-exercises resource, see Cazes and al. 2005). We consider these resources as artefacts, likely to become instruments for the teacher through an instrumental genesis process. Our aim is to observe and analyse this process. EEB can also become instruments for the students (Gueudet 2006); we do not study that aspect here. Applying the instrumental approach to teachers requires a specific theoretical work. We retained as a starting point the notion of scenario. We present in part 2 these theoretical choices. In part 3 we expose methodological elements: the data collected, the main EEB used by the teachers and the description grid of the scenarios. Part 4 is devoted to our main results and first observations of scenarios evolutions, interpreted as part of the genesis process.

## THEORETICAL FRAME: INSTRUMENTAL GENESIS AND SCENARIOS

Several research works study the complex practices of teachers with ICT (Monaghan 2004). Trouche (2004) introduces the notion of instrumental orchestration, which is associated with the instrumental approach and takes into account the teacher's action with ICT. But the orchestration aims at analysing how the teacher manages the students' instrumental genesis; it does not describe the teacher's own genesis.

Analysing teachers' instrumental genesis requires a specific theoretical reflection. It means to describe schemes, or techniques, considered as the visible part of the schemes developed by teachers (Guin and Trouche 2002). A first step is thus to observe if the teachers develop stable behaviours with the EEB, likely to be organised by an underlying scheme.

So we need to answer the following question: what does it mean to observe and describe teachers' practices with an EEB? For this purpose, we used the notion of scenario. Trouche (2004) defines a scenario in use as the organisation of a situation in a given environment. The scenario encompasses the management of the situation itself, and the management of the artifacts (orchestration). We consider that an EEB requires a specific approach because it can intervene also in the management of the situation.

Therefore we turned to the notion of scenario studied by computer science specialists, notably Pernin and Lejeune (2004), whose work is inspired by the Instructional Management Systems Learning Design (IMS LD) (Koper 2001) concepts and methods. We needed to adapt their approach, because we aim at describing and analyzing teachers' practices and not at developing technology enhanced learning techniques for them. We now present the theoretical frame elaborated by Pernin to study scenarios and the way we use it.

For Pernin and Lejeune, a learning scenario « represents a description, made a priori or a posteriori, of the progress of a learning situation at a given level, or learning unit, whose goal is to ensure the appropriation of a precise set of knowledge. A scenario describes roles, activities and also knowledge resources, tools and services necessary to the implementation of each activity. » (Lejeune and Pernin 2004). In our study, scenarios' designers are the teachers and users are the students. Pernin and Lejeune propose a list of six criterions to elaborate a taxonomy of scenarios.

**Purpose:** a scenario can be **scheduled** (developed a priori in order to define a learning situation) or **descriptive** (describing what actually happened during the learning situation). In our study, the scenarios can be scheduled or descriptive.

**Granularity:** Pernin and Lejeune introduce three granularity levels: the level of an elementary activity, the level of a sequence or composite activity and eventually the level of a structuration unit. Teachers may plan to use ICT for a one hour session to pursue a precise mathematical goal or may construct a sequence (a set of sessions with a same mathematical content).

**Constraint:** a **constrained** scenario describes precisely the activities to be realized, leaving a reduced initiative to the actors of the learning situation. An **open** scenario let them organize their progression; the main lines of the activities to be performed are described. At last, an **adaptable** scenario is an open scenario which can be modified by the actors of the learning situation. We can analyze the teacher's scenarios with this criterion: are they constrained or are students free to thumb through the EEB? For us, a scenario can never be adaptable because a student doesn't have the possibility to modify a scenario established by a teacher.

P**ersonalization**: a scenario can be **generic** (its execution is always the same) or **adaptive** (taking into account the user's profile and allowing execution of personalized scenarios).We use this criterion to sort scenarios depending on whether teachers chose to differentiate their students a priori, for instance with personalized menus.

F**ormalization and reification:** we only consider here **informal scenario**s, empirically created by teachers (and not **formalized scenarios** using an educational modelling language or **computable scenarios** written in a computable educational modelling language) and **concrete or contextualized scenarios (**and not **abstract scenarios)**. Scenarios produced by teachers are always concrete, which means they correspond to a precise context.

The formalization and reification criteria have always the same value, while the others are variable criteria: they can vary during the scenario life or, for the same teacher, they can be different from a scenario to another. We recall these variable criteria in the following table.

**Table 1: our use of Pernin's variable criteria to characterize scenarios.**

| | | | |
|---|---|---|---|
| **Purpose** | Scheduled | Descriptive | |
| **Granularity** | Session | Sequence | Structuration unit |
| **Personalization** | Generic | Adaptive | |
| **Constraint** | Constraint | Open | ~~Adaptable~~ |

All these aspects must be taken into account in the descriptions. But a didactical study requires also additional features. We describe them in the next part.

## METHODOLOGY

The data collected in 2005-2006 have two main origins: the five teachers of our group, and twenty other teachers, who followed a training course about a specific EEB called Mathenpoche ("Maths in the Pocket", shortened as MEP). For the sake of brevity, we do not mention here the questionnaire results. Moreover, we are now planning interviews useful to complement the questionnaires answers.

The teachers of our group are three primary school teachers (students between 9 and 11 years old); and two secondary school mathematics teachers (students between 12 and 15). Some are experienced with EEBs and some discover it, we explain it further in the next part. They were chosen to obtain this variety of contexts and experiences. During the whole 2005-2006 year all of them described precisely their uses of EEB, including: their preparation work, the description of what they planned, the description of what actually happened, the work done afterwards. Some of their sessions were observed, but only a few of them; thus the results presented in part 4 stem from the analysis of their own descriptions, formulated along a grid elaborated by the group.

The main resource used by the teachers of the group was MEP. We do not intend here to promote or dismiss any resource or kind of resource. MEP was the teachers' choice, thus a short insight of MEP's main features, and the presentation of the grid are necessary to understand our analysis.

**MEP's main features**

MEP offers exercises covering grade 6 to 9 curriculum, and a small part of grades 5 and 10. They are organised in sets of five or ten exercises with a common structure or theme; a given set of five exercises is identified by its title. Within an exercise set, the screen proposed to the student displays the text of one of the exercises (called a "problem", or a "question") with an answer zone to be filled in; the mark of the student (out of five or ten at the end of a given set of exercises); a button "calculator" providing access to a simple calculator; and a button "help" providing access to a full, explained solution of a similar exercise (always the same for a given exercise set). The expected answer can be numerical; some exercises offer multiple choices. After submitting their answer, the students receive a "Right" or "Wrong" feedback from the computer. Moreover, one or two detailed solutions of the exercise are displayed if their answer is right, or after a second wrong answer.

MEP is a free resource; some teachers use it with their classes without being registered as MEP's users. But registering as MEP's user opens more possibilities. The registered teachers inscribe their students into MEP. Each student is identified by a login and a password. Then the teacher chooses the exercise sets (these sets cannot be broken into smaller pieces) s/he wants to present to the students. The choice can be the same for the whole class, or different for subgroups, or even for individuals. The path of the students among the exercise sets selected can be left free or can be restricted according to the teacher's preferences: for example the second set could be offered only after a given threshold mark has been reached for the first. The registered teacher can follow the students' activity directly, through a special screen, during the sessions, or later by reading the "sessions' sheet" (see figure 1 below). It provides for each student (or group of students working on the same computer) the exercise sets tackled during the session, the time spent, the average, maximum and minimum mark obtained, and for each exercise, the success (green or medium grey in a not coloured print) or failure (red/dark grey; the question not tackled are displayed in blue/light grey).

**Figure 1 Extract of a session's sheet.**

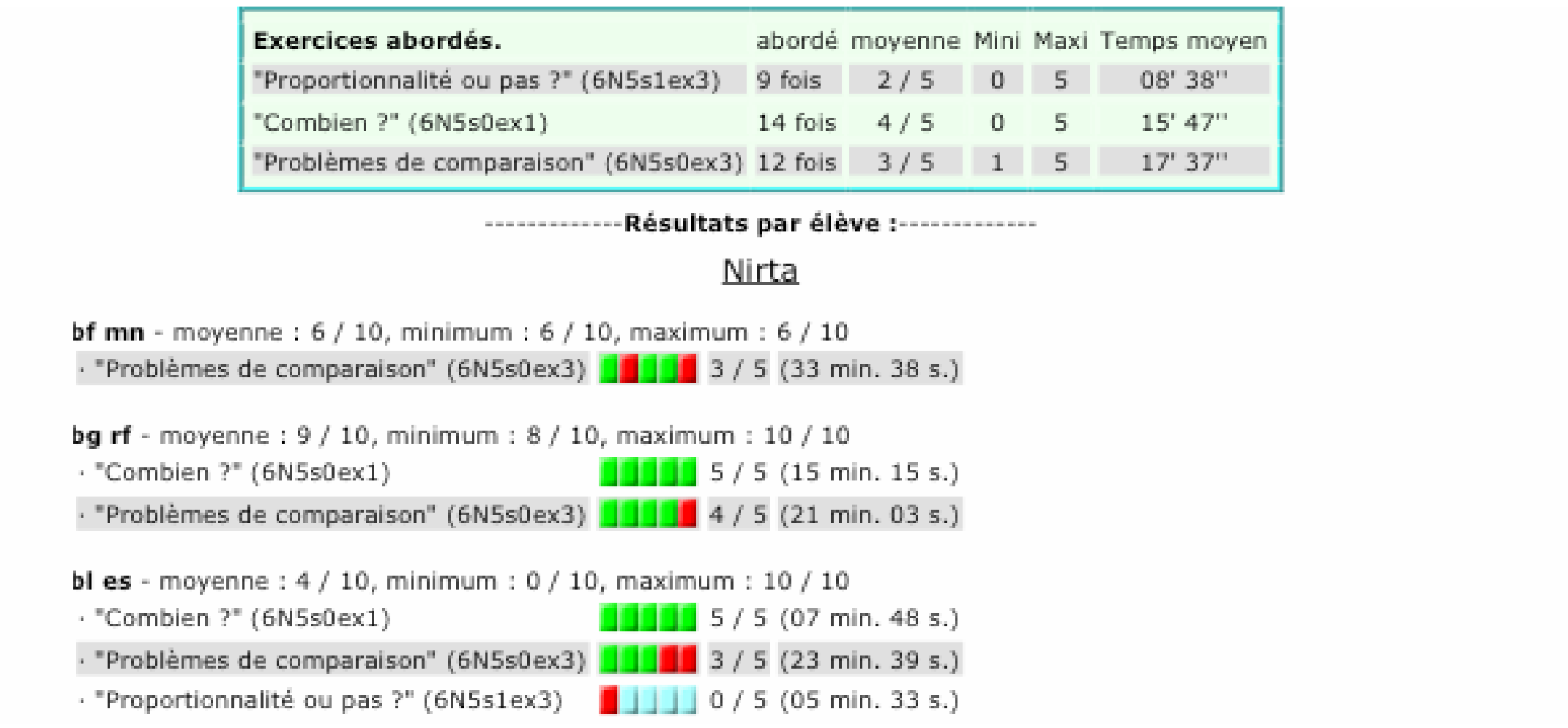

| Exercices abordés. | abordé | moyenne | Mini | Maxi | Temps moyen |
|---|---|---|---|---|---|
| "Proportionnalité ou pas ?" (6N5s1ex3) | 9 fois | 2 / 5 | 0 | 5 | 08' 38'' |
| "Combien ?" (6N5s0ex1) | 14 fois | 4 / 5 | 0 | 5 | 15' 47'' |
| "Problèmes de comparaison" (6N5s0ex3) | 12 fois | 3 / 5 | 1 | 5 | 17' 37'' |

-------------Résultats par élève :-------------

Nirta

**bf mn** - moyenne : 6 / 10, minimum : 6 / 10, maximum : 6 / 10
· "Problèmes de comparaison" (6N5s0ex3) 3 / 5 (33 min. 38 s.)

**bg rf** - moyenne : 9 / 10, minimum : 8 / 10, maximum : 10 / 10
· "Combien ?" (6N5s0ex1) 5 / 5 (15 min. 15 s.)
· "Problèmes de comparaison" (6N5s0ex3) 4 / 5 (21 min. 03 s.)

**bl es** - moyenne : 4 / 10, minimum : 0 / 10, maximum : 10 / 10
· "Combien ?" (6N5s0ex1) 5 / 5 (07 min. 48 s.)
· "Problèmes de comparaison" (6N5s0ex3) 3 / 5 (23 min. 39 s.)
· "Proportionnalité ou pas ?" (6N5s1ex3) 0 / 5 (05 min. 33 s.)

Naturally, the teachers in our group are registered as MEP's users; it is the main resource they used.

## A grid for the description of scenarios' features

The grid for the description of scenarios' features was developed with the teachers of the group. It includes Pernin and Lejeune's criteria, but also all the other features that appeared in the teachers' informal descriptions.

They systematically focus their descriptions on the mathematical knowledge to reach. The source of theirs reflections is the knowledge to be taught. The question of using an EEB and the way to do so comes next. It can lead teachers to build their scenarios either at the session or at the sequence granularity level. These two possibilities appear in the following grid.

**Table 2: grid for the description of scenarios' features**

| | |
|---|---|
| 1. Work plan | 1.1 Resources: chosen EEB, other ICT, other supports |
| | 1.2 For a sequence : sessions distribution and articulation, number, time |
| | 1.3 Mathematical content, objectives |
| | 1.4 Session(s) type(s) : discovery, remedial, evaluation… ; precise the EEB function |
| | 1.5 For a sequence: references to the EEB in paper and pencil sessions |
| 2.Teacher's interventions during computer session(s) | 2.1 Content |
| | 2.2 Support |
| 3. Students' activities | 3.1 Written notes expected during computer session (s) and written documents provided |
| | 3.2 Alone / pair / group work during computer session(s) |
| | 3.3 Work with computer out of class |
| | *3.4 Differentiation (Personalization)* |
| | *3.5 Progression: imposed, with thresholds, free (Constraint)* |

The grid comprises components stemming from our didactical approach and others corresponding to Pernin's variable criteria (italicized in the table). The teachers filled in such a grid for each of their uses of an EEB. We now present the main results, in terms of scenario's evolution observed.

## GENESIS PROCESSES OBSERVED IN THE GROUP

We present in the following table basic numerical data about the use of web resources by the members of our group.

**Table 3: Use of web resources during the year 2005-2006**

| | Teacher 1 (Secondary school) | Teacher 2 (Secondary school) | Teacher 3 (Primary school) | Teacher 4 (Primary school) | Teacher 5 (Primary school) |
|---|---|---|---|---|---|
| Number of scenarios | 7 | 6 | 5 | 6 | 4 |
| Use of the resource in class (hours) | 15 | 10 | 15 | 9,5 | 8 |

These data help to figure out what happened in the classes; however, we are more interested in the evolutions observed. Two teachers in the group were familiar with the use of MEP, it was indeed their third year using MEP in their classes. They both used MEP, but also another EEB completely new for them. Two teachers already met MEP at the end of the preceding year (and used no other resource); we consider that they started being familiar with it around January, when they overcame all the technical difficulties. One teacher discovered MEP by joining the group (she also used two other resources). According to the description of scenarios provided by the teachers, we retain several kinds of noticeable evolutions during the year 2005-2006.

From isolated sessions to sequences

Two main granularity levels appear in the teachers' descriptions: the session and the sequence. Less frequently, other time scales are mentioned, corresponding to specific uses: training out of class, programmed for several weeks or even months; or individual help. But the work intended for the whole class is thought of on the session, or on the sequence granularity.

The ability to plan at the sequence's scale appears as a result of the genesis process. For unfamiliar resources, 12 scenarios are proposed; 6 are planned for one-hour sessions, and 6 for sequences. For familiar resources, 16 scenarios are proposed; 13 for sequences and only 3 for one-hour sessions. The discussions in the group confirm these observations. With a new resource, the teachers have a wide mathematical content to discover, they can not rely on well-known exercises to build a progression. For example, when the two secondary school teachers discovered the resource WIMS (Web International Multipurpose Server, http://wims.unice.fr), they retained a ready-made session of training about powers of ten. Discovering all the exercises in this theme to build a progression was a too tough work (understanding the resource features already required two hours). Thus they only organized this isolated training session; it was related to their traditional teaching on the theme, but not really integrated in it.

On the opposite, a familiar resource is an instrument for the teacher just like the textbook for example. For a given mathematical content, the teacher plans a sequence, and foresees where the EEB will intervene if he or she already knows that the EEB proposes interesting exercises about this content. For example, for the theme "parallel and perpendicular straight lines" for grade 6 pupils, a sequence of 11 hours was planned, including 2 hours on MEP, 2 hours on a dynamic geometry software, and 7 paper and pencil hours. All these kinds of sessions were strongly intertwined; the aim of the work on MEP was described as "observing figures and discovering first proof principles", and a synthesis was organized in class after it.

A wide range of scenarios

The scenarios elaborated by the teachers of the group display an important variety of features; we emphasise here three main directions where evolutions were observed during the year.

About the function of the EEB: the teachers unfamiliar with an EEB mostly use it for training on technical abilities, like drill and practice software. Another use planned by unfamiliar teachers is the projection of the solving of an exercise with a video for the whole class. Both cases share a common feature. They avoid the discussion in class about exercises done by the students on the computer. The teachers in the group unfamiliar with MEP expressed it several times: they feared to organise a discussion after the work on MEP, because the students met different exercises on the computer, and the teacher found difficult to discuss in class contents that some students perhaps had not met. This obstacle was overcome mostly thanks to the session's sheet, that allows the teacher to identify the exercises tackled by all the students. Then it became possible to use MEP to discover new techniques and properties; the session's sheet also made possible the use of MEP for evaluation purposes. The evaluation function does not appear in the scenarios organised by teachers unfamiliar with MEP.

About the written notes expected during the computer sessions: written notes seem to be expected more frequently with familiar resources. It is naturally connected with the preceding point: written notes do not seem so necessary during a drill and practice session, because there will be no further work on such a session. Moreover, it is easier to choose the appropriate written note for familiar exercises. It is indeed difficult for the students to write everything they do on the computer. Thus a choice is necessary. For example, at the end of the year, one teacher asked the students to write down the full solution of one MEP's exercise for each exercise set (five exercises) tackled.

About general organisation choices: the whole class can work at the same moment on the computer, students working by pairs. Or subgroups can work individually on the computer, while the other students have a paper and pencil work. Out of class work has also been planned by one primary school teacher. These choices are naturally strongly related with technical conditions. A half class on the computer, and the other half class on paper is only possible if the computers are situated in a room with enough additional tables. And the work out of class was only programmed by teachers who know their pupils have an Internet access at home, or at a nearby library. But the increasing confidence of the teacher with the EEB seems also crucial. A good knowledge of the exercises permits for example to retain exercises interesting for exchanges within a pair (for example exercises with two possible solutions).

The discussions with other teachers during the group meetings also played an important part in the development of various scenarios. The only teacher who was not able, for geographical reasons, to join the group's meetings always proposed similar scenarios (sequences of around ten sessions, with one or two hours on MEP for drill on technical exercises without written notes). On the opposite, the other teachers clearly influenced each other; two teachers working in the same secondary school regularly worked together to prepare their scenarios, and even ended up with a scenario associating two of their classes split in three level groups.

### More personalization and constraints

The scenarios planned display more and more personalization and constraints. Two ways of personalization and constraints are used: one is the use of the possibilities offered by the resource (mostly MEP); the other is the formulation of advice by the teacher.

Some teachers program a different content for different groups of students, even sometimes for individuals. This is the personalization explicitly intended by the designers of the resource. At the end of the year, three of the five teachers chose on the opposite to propose a wide range of exercises, the same for all students. The personalization was then done through individual advice on the exercises to tackle.

Similarly about constraints, some teachers use the constraints proposed by the resource: they program sessions were exercise 2 can only be tackled after exercise 1 (imposing a threshold mark was quickly let down). But at the end of the year, they seem to prefer proposing all the exercises together, and formulating their own constraints, the most frequent being about the number of attempts on the same exercise: never more than two attempts, whatever the mark is.

In both cases, the teachers first become acquainted with the possibilities offered by the resource; then they elaborate their own possibilities, sometimes not expected by the designers of the EEB. They create their own instrument from the EEB.

## CONCLUSION AND PERSPECTIVES

We clearly observed evolutions of the teachers' practices. These evolutions are linked with their mastery of the technical features of the EEB, but mostly with their knowledge of its mathematical content.

We demonstrated that the practices of the group's novices teachers evolved towards the experts' ones. It allows us to consider the experts' practices as stable behaviours. The teachers of the group built an instrument from the EEBs they used.

It supports the interpretation of the practices' evolutions as teachers' instrumental geneses. Thus the next step of our work is to describe techniques considered as the visible part of the schemes developed by teachers. First, we have to identify marks of instrumentation and instrumentalization phenomena in teachers' practices with EEBs.

Tasks and techniques in the teacher's action have been studied by Sensevy and al. (2005). These authors identify tasks and techniques in the teachers' practices through direct observation and through the teachers' descriptions of their own practices. In a similar way, we are now trying to identify, in the scenarios described by the teachers but mostly through class observation, techniques instrumented by the EEB. For example, the EEB seems to provide instrumented techniques for the task "managing the students heterogeneity", because it permits for example 1. an instrumented diagnosis of students' knowledge by choosing appropriate exercises, 2. an instrumented choice of contents by programming individual menus and 3. an instrumented heterogeneity's management, allowing for example teachers to follow the students' individual work.

But studying the instrumented techniques requires an analysis focused on a particular mathematical content, allowing a comparison, on more precise tasks, of the teachers' paper and pencil techniques and their techniques with the EEB.